\documentclass[12pt]{article}
\usepackage{a4}
\usepackage{verbatim}
\usepackage{amsmath}
\usepackage{amssymb}
\usepackage{amsthm}

\usepackage{definitions}
\usepackage[dvips]{graphicx}

\author{Jaros\l{}aw Buczy\'nski\footnote{
The article is a part of the research project N20103331/2715 
funded by Polish financial means for science in years 2006-2008. 
Author's e-mail: jabu@mimuw.edu.pl}
}

\title{Toric Legendrian subvarieties}
\date{March 30, 2007}

\begin{document}
\maketitle

\begin{abstract}
We give the full classification of smooth toric Legendrian subvarieties in $\P^{2n-1}$.
We also prove that under some minor assumptions the group of linear automorphisms preserving a given 
Legendrian subvariety preserves the contact structure of the ambient projective space.
\end{abstract}

\section{Introduction}
We are interested in Legendrian subvarieties of projective space $\P^{2n-1}$
and here we recall the definition:

\begin{defin}
Let $\omega$ be a symplectic form on $V=\C^{2n}$.
A subvariety $X \subset \P(V)$ is \emph{Legendrian}, 
if for each smooth point of its affine cone $\hat{X}$ the tangent space to $\hat{X} \subset V$ at this point is
Lagrangian, i.e.~maximal isotropic with respect to $\omega$. 
\end{defin}

The importance and relation of Legendrian varieties 
to the problem of classifying contact Fano varieties is briefly
explained in \cite[\S2]{jabu} and the reference therein. 
There is another equivalent definition: 
a subvariety is Legendrian if its tangent bundle is contained in the contact distribution 
on $\P^{2n-1}$. 
It is explained and presented for example in \cite{landsbergmanivel04}.

\begin{defin}\label{definition_decomposable}
Let $V_1$ and $V_2$ be two symplectic vector spaces and 
let $X_1\subset \P(V_1)$ and $X_2\subset \P(V_2)$ be two Legendrian subvarieties.
Now assume $V:=V_1 \osum V_2$ and $X:=X_1 * X_2 \subset \P(V)$, i.e.~$X$ is the joint of $X_1$ and $X_2$ 
meaning the union of all lines from $X_1$ to $X_2$. 
Now, clearly, the affine cone of $X$ is the product $\hat{X_1} \times \hat{X_2}$ 
(where $\hat{X_i}$ is the affine cone of $X_i$). 
In such a case we say that $X$ is a \emph{decomposable Legendrian variety}. 
We say that a Legendrian subvariety in $V$ is \emph{indecomposable} 
if it is not of that form for any non-trivial symplectic decomposition $V = V_1 \osum V_2$.
\end{defin}

We state the following conjecture:

\begin{con}\label{theorem_automorphisms_are_symplectic}
Let $X\subset \P^{2n-1}$ be an irreducible indecomposable 
Legendrian subvariety and let $G< \P\Gl_{2n}$ be a connected subgroup preserving $X$. 
Then $G$ is contained in the image of the natural map $\Sp_{2n} \ra \P\Gl_{2n}$.
\end{con}

We prove the conjecture in several cases.

It is quite natural to believe, that if a linear map preserves a form on a big number of linear
subspaces, then it actually preserves the form (at least up to scalar). 
With this approach,
\cite[cor. 6.4]{janeczkojelonek} proved the conjecture in the case
 where the image of $X$ under the Gauss map is 
non degenerate in the Grassmannian of Lagrangian subspaces in $\C^{2n}$.
Unfortunately, this is not enough -  for example $\P^1 \times Q_1 \subset \P^5$ 
has a degenerate image under the Gauss map 
and this is one of the simplest examples of smooth Legendrian subvarieties.

In section \ref{section_proof_of_the_conjecture} we prove:

\begin{theo}\label{corollary_that_conjecture_is_true_for_smooth}
If $X\subset \P^{2n-1}$ is a smooth Legendrian subvariety which is not a linear subspace
and $G< \P\Gl_{2n}$ is a connected subgroup preserving $X$,
then $G$ is contained in the image of the natural map $\Sp_{2n} \ra \P\Gl_{2n}$.
\end{theo}

The theorem is applied to classify smooth toric Legendrian subvarieties. 
So in section \ref{section_toric} we choose appropriate coordinates to 
reduce this problem to some combinatorics (for surface case --- see section \ref{section_toric_surfaces})
and some elementary geometry of convex bodies (for higher dimensions --- 
see section \ref{section_higher_dimensional}). 
Eventually we get:

\begin{theo}
Every smooth toric Legendrian subvariety in a projective space is isomorphic 
to one of the following:
\begin{itemize}
\item a linear subspace,
\item  $\P^1 \times Q_1 \subset \P^5$,
\item  $\P^1 \times Q_2 \simeq \P^1 \times \P^1 \times \P^1 \subset \P^7$
\item  or $\P^2$ blown up in three non-colinear points.
\end{itemize}
\end{theo}
For proofs see corollaries \ref{smooth_toric_surfaces} and \ref{smooth_toric_varieties}.
The linear subspace is not really interesting, 
the products $\P^1 \times Q_1$ and $\P^1 \times Q_2$  are well known 
(see for example \cite{landsbergmanivel}, \cite{mukai}, \cite{jabu}).
The last case of blow up is not yet described in the literature.

\subsection{Acknowledgements}
The author would like to thank especially Jaros\l{}aw Wi\'s{}niewski for his comments and support.
Also I am grateful to Micha\l{} Kapustka, Grzegorz Kapustka, Adrian Langer, Zbigniew Jelonek, Michel Brion 
and the anonymous referees for their precious remarks and hints.

\section{Projective automorphisms of a Legendrian \\ variety}\label{section_automorphisms}

We would like to make sure that the torus acting on a Legendrian variety acts 
via symplectic automorphisms of the associated symplectic vector space. 
Therefore in this section we partially  answer the question,
whether a projective automorphism of a Legendrian variety must be symplectic.

One obvious remark is that homotheties act trivially on $\P(V)$, but in general are not symplectic.
Therefore, it is more convenient to think of conformal symplectomorphisms:

\begin{defin}
A linear automorphism $\psi$ of a symplectic vector space $(V,\omega)$ is called
a \emph{conformal symplectomorphism}
if $\psi^* \omega = c \omega$ for some constant $c \in \C^*$. 
We denote by $\cSp (V)$ the group of all conformal symplectomorphisms of $V$ and
by $\gotcsp(V)$ the tangent Lie algebra.
\end{defin}

Consider the following example:
 
\begin{ex}\label{example_degenerate}
Let $X\subset \P(V)$ be a Legendrian subvariety contained in a hyperplane. 
Then by \cite[thm. 3.4]{jabu} the vector space $V$ admits a symplectic decomposition 
$V= W \oplus H$, where $\dim H =2$ and $X$ is a cone over $X \cap \P(W)$ with a vertex in $\P(H)$.
Then a linear automorphism that acts as $\lambda_1 \Id$ on $W$ and as $\lambda_2 \Id$ on $H$
preserves $X$ and is definitely not a conformal symplectomorphism, unless $\lambda_1 = \pm \lambda_2$.
\end{ex}

Therefore it is clear, that if we hope for a positive answer to the question posed in the beginning of
this section we must assume that our Legendrian variety is non degenerate (i.e.~not contained in 
any hyperplane). 
Another natural assumption is that $X$ is irreducible 
--- one can also easily produce a counterexample if we skip this assumption.
Yet still this is not enough.

Let $X= X_1* X_2\subset \P(V_1 \oplus V_2)$ be a decomposable Legendrian variety. 
Then as in example \ref{example_degenerate} we can act via $\lambda_1 \Id_{V_1}$ on $V_1$ and 
via $\lambda_2 \Id_{V_2}$ on $V_2$ - such an action will preserve $X$ 
and again in general it is not conformal symplectic. 
This explains that the assumptions of our conjecture \ref{theorem_automorphisms_are_symplectic}
are necessary.

Yet we must note, that decomposable varieties get very singular, unless they are linear space.
Our main interest is in smooth Legendrian subvarieties.



\subsection{Decomposable varieties}

We prove an easy proposition.

\begin{prop}\label{proposition_decomposable}
Let $S_1$ and $S_2$ be two smooth algebraic varieties and let $X \subset S_1 \times S_2$ 
be a closed irreducible subvariety.
Let $X_i \subset S_i$ be the closure of the image of $X$ under the projection $\pi_i$ onto $S_i$. 
Assume that for a Zariski open dense subset of smooth points $U\subset X$
we have that the tangent bundle to $X$ decomposes as
$T X|_U = (T X \cap \pi_1^* T {S_1})|_U \oplus (T X \cap \pi_2^* T {S_2})|_U$ a sum of two vector bundles.
Then $X= X_1 \times X_2$.
\end{prop}

\begin{prf}
Since $X$ is irreducible, so is $X_1$ and $X_2$ and clearly $X \subset X_1 \times X_2$.
 So it is enough to prove that $\dim X_1 + \dim X_2 = \dim X = \dim U$.
But the maps $\ud (\pi_i|_U)$ are surjective onto $T X \cap \pi_i^*T {S_i}$ and hence by \cite[thm III.~10.~6]{hartshorne}
$$
\dim X_1 + \dim X_2 = \rk (T X \cap \pi_1^* T S_1)|_U + \rk (T X \cap \pi_2^* T{S_2})|_U = \rk T X|_U = \dim X.
$$
\end{prf}


\subsection{Weks-symplectic matrices}

Fix a basis of $V$ and recall
 that a matrix $g\in \gotgl(V)$ is in the symplectic algebra $\gotsp(V)$ if and only if 
$$g^TJ + Jg =0$$
where $J$ is the matrix of the symplectic form in the given basis. 
We want to define a complementary linear subspace to $\gotsp(V)$:
\begin{defin}
A matrix $g\in \gotgl(V)$ is \emph{weks-symplectic}\footnote{
A better name would be \emph{skew-symplectic} or \emph{anti-symplectic}, but these are reserved for some different notions.
}
if and only if it satisfies the equation:
$$
g^TJ - Jg =0.
$$
The vector space of all weks-symplectic matrices will be denoted by $\gotasp(V)$
(even though it is not a Lie subalgebra of $\gotgl(V)$).
\end{defin}

We immediately see that a matrix is weks-symplectic if and only if it corresponds to 
a linear endomorphism $g$, 
such that for every $u,v\in V$:
\begin{equation}\label{weks_symplectic_equation}
\omega(gu, v) - \omega (u, gv) =0.
\end{equation}
This is a coordinate free way to describe $\gotasp(V)$.

\medskip

From now on we assume that our basis is symplectic, which means that the matrix of the symplectic form
is of the following block form:
$$J= \left( \begin{array}{cc}
0&\Id_n\\
-\Id_n&0
\end{array} \right) .$$
In particular $J^2 = -\Id_{2n}$.

\begin{rem}\label{weks_corresponds_to_skew}
For a matrix $g\in \gotgl(V)$ we have:
\begin{itemize}
\item[(a)] $g \in \gotsp(V) \  \iff  \  Jg$ is a symmetric matrix;
\item[(b)] $g \in \gotasp(V) \  \iff  \  Jg$ is a skew-symmetric matrix.
\end{itemize}
\end{rem}
\noprf

Note that if $g\in \gotgl(V)$, then we can write:
$$
g = \half( g + J g^T J) + \half (g - J g^T J)
$$
and the first component $g_+:= \half( g + J g^T J)$ is in $\gotsp(V)$, 
while the second $g_- := \half (g - J g^T J)$ is in $\gotasp(V)$. 
Obviously, this decomposition corresponds to expressing the matrix $Jg$ 
as a sum of symmetric and skew-symmetric matrices.

We list some properties of $\gotasp(V)$:
\begin{prop}\label{properties_of_asp}
Let $g, h \in \gotasp(V)$. The following properties are satisfied: 
\begin{itemize}
\item[(i)] Write the additive Jordan decomposition for $g$:
$$
g = g_s + g_n 
$$
where $g_s$ is semisimple and $g_n$ is nilpotent.
Then both $g_s \in \gotasp(V)$ and $g_n \in \gotasp(V)$.

\item[(ii)] For $\lambda \in \C$, 
denote by $V_{\lambda}$ the \mbox{$\lambda$-eigenspace} of $g$.
For any $\lambda_1, \lambda_2 \in \C$ two different eigenvalues  
$V_{\lambda_1}$ is $\omega$-perpen\-di\-cu\-lar to $V_{\lambda_2}$.

\item[(iii)] If $g$ is semisimple, then each space $V_{\lambda}$ is symplectic,
i.e.~the form $\omega|_{V_{\lambda}}$ is non degenerate.

\end{itemize}
\end{prop}
\noprf

\subsection{About conjecture \ref{theorem_automorphisms_are_symplectic} 
and the proof of theorem \ref{corollary_that_conjecture_is_true_for_smooth}}
\label{section_proof_of_the_conjecture}

Let $X'\subset \P(V)$ be an irreducible, indecomposable Legendrian subvariety, 
let \mbox{$X$ be} the affine cone over $X'$  and $X_0$ be the smooth locus of $X$. 
Assume that $G$ is the maximal connected subgroup in $\Gl_{2n}$
preserving $X$. Let $\gotg< \gotgl_{2n}$ be the Lie algebra tangent to  $G$.
To prove the conjecture it would be enough to show that $\gotg$ is contained in the Lie algebra $\gotcsp_{2n}$
tangent to conformal symplectomorphisms, i.e.~the Lie algebra spanned by $\gotsp_{2n}$ and the identity
matrix $\Id_{2n}$.

\begin{theo}\label{theorem_evidence_for_conjecture}
With the above notation the following properties hold:
\begin{itemize}
\item[I.]
The underlying vector space of $\gotg$ decomposes into symplectic and weks-symplectic part:
$$
\gotg = \big(\gotg \cap \gotsp(V)\big)  \oplus \big(\gotg \cap \gotasp(V)\big). 
$$ 
\item[II.]
If $g \in \gotg \cap \gotasp(V)$, then $g$ preserves every tangent space to $X$:
$$
\forall {x\in X_0} \quad g(T_x X) \subset T_x X
$$
and hence also
$$
\forall t\in \C \quad \forall x\in X_0 \qquad T_{\exp(tg)(x)} X = \exp(tg)(T_x X) = T_x X.
$$
\item[III.]
If $g \in \gotg \cap \gotasp(V)$ is semisimple,
then $g=\lambda \Id$ for some $\lambda \in \C$.
\item[IV.]
Assume $0 \ne g \in \gotg \cap \gotasp(V)$ is nilpotent and let $m \ge 1$ be an integer such that $g^{m+1} = 0$
and $g^{m}\ne 0$.
Then $g^m(X)$ is always non-zero and
is contained in the singular locus of $X$.
In particular, $X'$ is singular.
\end{itemize}
\end{theo}

In what follows we prove the four parts of theorem \ref{theorem_evidence_for_conjecture}.

\subsubsection{I.~Decomposition into symplectic and weks-symplectic part}
\begin{prf}
Take $g \in \gotg$ to be an arbitrary element. 
Then for every $x \in X_0$  one has 
$$g (x) \in T_x X$$
and therefore:
$$
0 = \omega \big( g(x), x \big) = x^T g^T J x  = \half x^T \left( g^T J - J g \right) x
$$

Hence the quadratic polynomial  $f(x):=x^T ( g^T J - J g) x$ is identically zero on $X$ and hence it is in the ideal of $X$.
Therefore by maximality of $G$ and \cite[cor. 5.5]{jabu} 
the map $J \left(g^T J - J g\right)$ is also in $\gotg$. But
$$
J \left(g^T J - J g\right) = J g^T J + g,
$$
so $Jg^TJ \in \gotg$ and both symplectic and weks-symplectic components $g_+$ and $g_-$ are in $\gotg$.
\end{prf}

From the point of view of the conjecture, the symplectic part is fine.
We would only need to prove that $g_- = \lambda \Id$.
So from now on we assume $g = g_- \in \gotasp(V)$.

\subsubsection{II.~Action on tangent space}
\begin{prf}
Let $\gamma_t:=\exp(t g)$ for $t \in \C$.
Then $\gamma_t \in G$ and hence it acts on $X$.
Choose a point $x\in X_0$ and two tangent vectors in the same tangent space $u,v \in T_x X$.
Then clearly also $\gamma_t(u)$ and $\gamma_t(v)$ are contained in one tangent space, 
namely $T_{\gamma_t(x)} X$. Hence:
$$
0= \omega \left(\gamma_t(u),\gamma_t(v)\right) = 
\omega \big(
(\Id_{2n} + tg  +\ldots)u, 
(\Id_{2n} + tg  +\ldots)v
\big) = 
$$\nopagebreak
$$
= \omega(u,v) + t \big( \omega(gu, v) + \omega(u, gv)\big) + t^2(\ldots).
$$
In particular the part of  the expression linear in $t$ vanishes, hence $\omega(gu, v) + \omega(u, gv) =0$.
Combining this with equation (\ref{weks_symplectic_equation}) we get that:
$$
\omega(gu, v) = \omega(u, gv) = 0
$$
But this implies that $gu \in (T_x X)^{\perp_{\omega}} = T_x X$.
Therefore $g$ preserves the tangent space at every smooth point of $X$ 
and hence also $\gamma_t$ does preserve that space.
\end{prf}

\subsubsection{III.~Semisimple part}
Since $G$ is an algebraic subgroup in $\Gl(V)$, 
then $\gotg$ has the natural Jordan decomposition inherited from $\gotgl(V)$, 
i.e.~if we write the Jordan decomposition for $g= g_s +g_n$, then $g_s, g_n \in \gotg$ 
(see \cite[thm 15.3(b)]{humphreys2}).
Therefore by proposition \ref{properties_of_asp} (i),
it would be enough, to establish the conjecture, to prove that
for $g\in \gotg \cap \gotasp(V)$ we have $g_s = \lambda \Id_{2n}$ and $g_n=0$.

Here we deal with the semisimple part.

\begin{prf}
Argue by contradiction.
Let $V_1$ be an arbitrary eigenspace of $g$ and let $V_2$ be the sum of the other eigenspaces. 
If $g \ne \lambda \Id_{2n}$, then both $V_1$ and $V_2$ are non-zero and by 
proposition \ref{properties_of_asp} (ii) and (iii) they are $\omega$-perpendicular, 
complementary symplectic subspaces of $V$. Let $x\in X_0$ be any point. 
Since $g$ preserves $T_x X$ by part II.~it follows
that  $T_x X= (T_x X\cap V_1) \osum (T_x X \cap V_2)$. 
But then both $(T_x X\cap V_i) \subset V_i$ are Lagrangian subspaces, hence have constant 
(independent of $x$) dimensions.
Hence $T_x X_0 = (T_x X_0 \cap V_1) \osum (T_x X_0 \cap V_2)$ is a sum of two vector bundles 
and  from proposition \ref{proposition_decomposable} we get that 
$X$ is a product of two Lagrangian subvarieties contradicting our assumption that $X'$ is indecomposable.
\end{prf}


\subsubsection{IV.~Nilpotent part --- $X'$ is singular}

\begin{lem}\label{limit_of_exptgv}
Assume $X' \subset \P(V)$ is any closed subvariety preserved by the action of $\exp(tg)$
for some nilpotent endomorphism $g\in \gotgl(V)$.
If $v$ is a point of the affine cone over $X'$ and $m$ is an integer such
that $g^{m+1}(v)=0$ and $g^{m}(v) \ne 0$, then the class of $g^{m}(v)$ is in $X'$.
\end{lem}

\begin{prf}
The class of $g^{m}(v)$ in the projective space $\P(V)$ is just the limit of classes of $\exp(tg)(v)$ as $t$ 
goes to $\infty$.
\end{prf}

\begin{lem}\label{smooth_point_implies_linear}
Assume $g\in \gotgl(V)$ is nilpotent and
 $g^{m+1}=0$, $g^{m} \ne 0$ for an integer $m\ge 1$.
Let $X\subset V$ be an affine cone over some irreducible projective subvariety in $\P(V)$,
which is preserved by the action of $\exp(tg)$, but is not contained in the set of the fixed points. 
Assume that this action preserves the tangent space $T_x X$ at every smooth point $x$ of $X$.
If there exists a non-zero vector in $V$ which is a smooth point of $X$ contained in $g^{m}(X)$,
then $X$ is a linear subspace.
\end{lem}

\begin{prf}
\emph{Step 0 - notation.}
We let $Y$ to be the closure of $g^m(X)$, so in particular $Y$ is irreducible.
By lemma \ref{limit_of_exptgv},
we know that $Y \subset X$.
Let $y$ be a general point of $Y$.  
Then by our assumptions $y$ is a smooth point of both $X$ and $Y$.

Next denote by 
$$
W_y:= (g^m)^{-1} (\C^* y).
$$  
You can think of $W_y$ as union of those lines in $V$
(or points in the projective space $\P(V)$), 
which under the action of $\exp(tg)$
converge 
to the line spanned by $y$ (or [y])\footnote{
This statement is not perfectly precise, though it is true on an open dense subset.
There are some other lines, which converge to $[y]$, namely those generated by $v \in \ker g^m$, 
but $g^k(v) = \lambda y$ for some $k<m$.
We are not interested in those points.
}
as $t$ goes to $\infty$ .
We also note that the closure $\overline{W_y}$ is a linear subspace spanned by an arbitrary element 
$v \in W_y$ and $\ker g^m$.

Also we let $F_y:= W_y \cap X$, so that:
$$
F_y:= (g^m|_X)^{-1} (\C^* y).
$$
Finally, $v$ from now on will always denote an arbitrary point of $F_y$.

\medskip

\emph{Step 1 - tangent space to $X$ at points of $F_y$.}
Since $y$ is a smooth point of $X$ also $F_y$ consists of smooth points of $X$. 
This is because the set of singular points is closed and $\exp(tg)$ invariant. 
By our assumptions  $\exp(tg)$ preserves every tangent space to $X$ 
and thus for every $v \in F_y$ we have:
$$
T_v X = T_{\frac{1}{t^m}\exp(tg)(v)} X = T_{\lim_{t\ra \infty} \left(\frac{1}{t^m}\exp(tg)(v)\right)} = T_y X.
$$
So the tangent space to $X$ is constant over the $F_y$ and in particular $F_y \subset T_y X$.

\medskip

\emph{Step 2 - dimensions of $Y$ and  $F_y$.}
From the definitions of $Y$ and $y$ and by step 1 we get that for any point $v\in F_y$:
$$
T_y Y \ = \ \im (g^m|_{T_v X}) =   \im (g^m|_{T_y X}).
$$
Hence $\dim Y = \dim T_y Y = \rk (g^m|_{T_y X})$.

Since $y$ was a general point of $Y$, we have that:
$$
\dim Y + \dim F_y = \dim X +1.
$$
So $\dim F_y = \dim \ker(g^m|_{T_y X}) +1$.

\medskip

\emph{Step 3 - the closure of $F_y$ is a linear subspace.}
From the definition of $F_y$ and step 1
we know that $F_y \subset T_y X \cap W_y$ and 
$$
 T_y X \cap \overline{W_y} =
 T_y X \cap \sspan\{v, \ker g^m \} = 
\sspan\{v, \ker (g^m|_{T_y X}) \}
$$
hence $\dim F_y = \dim T_y X \cap W_y$,
so the closure of $F_y$ is exactly $ T_y X \cap \overline{W_y}$ and clearly this closure is contained in $X$.
In particular $\ker (g^m|_{T_y X}) \subset X$. 

\medskip

\emph{Step 4 - $Y$ is contained in $\ker (g^m|_{T_y X})$.}
Let $Z$ be $X \cap \ker g^m$.
By step 3 we know that $\ker (g^m|_{T_y X}) \subset Z$. 
Now we calculate the local dimension of $Z$ at $y$:
$$
\dim \ker (g^m|_{T_y X}) \le
\dim_y Z \le \dim T_y Z \le 
\dim (T_y X \cap \ker g^m) = \dim \ker (g^m|_{T_y X}). 
$$
Since the first and the last entries are identical, we must  have all equalities. 
In particular the local dimension of $Z$ at $y$ is equal to the dimension of the tangent space to $Z$ at $y$.
So $y$ is a smooth point of $Z$ and therefore there is a unique component of $Z$ passing through $y$, 
namely the linear space $\ker (g^m|_{T_y X})$. 
Since $Y$ is contained in $Z$ (because $\im g^m \subset \ker g^m$) and $y \in Y$, 
we must have $Y \subset  \ker (g^m|_{T_y X})$.

\medskip

\emph{Step 5 - vary $y$.}
Recall, that by step 1 the tangent space to $X$ is the same all over $F_y$. 
So also it is the same on every smooth point of $X$, 
which falls into the closure of $F_y$. 
But by step 4, $Y$ is a subset of $ \ker (g^m|_{T_y X})$, 
which is in the closure of $F_y$ by step 3. 
So the tangent space to $X$ is the same for an open subset of points in $Y$.
Now apply again step 1 for different $y$'s in this open subset 
and we get that $X$ has constant tangent space on a dense open subset of $X$.
This is possible if and only if $X$ is a linear subspace, which completes the proof of the lemma.
\end{prf}

Now part IV.~of the theorem follows easily:

\begin{prf}
By the assumptions of the theorem $X$ is  not contained in any hyperplane, 
so in particular $X$ is not contained in $\ker g^m$. 
So by lemma \ref{limit_of_exptgv} the image $g^m(X)$ contains other points than $0$.
Next by lemma \ref{smooth_point_implies_linear} and part II.~of the theorem, since $X$ cannot be a linear subspace,
there can be no smooth points of $X$ in  $g^m(X)$.  
\end{prf}

\subsubsection{Smooth case}
We conclude that parts I., III.~and IV.~of theorem \ref{theorem_evidence_for_conjecture}
together with proposition  \ref{properties_of_asp} (i) and \cite[thm. 15.3(b)]{humphreys2}
imply theorem \ref{corollary_that_conjecture_is_true_for_smooth}. 
We only note that a smooth Legendrian subvariety is either a linear subspace or it is indecomposable.



\subsection{Some comments}\label{section_comments}

Conjecture \ref{theorem_automorphisms_are_symplectic} is now reduced  
to the following special case not covered by theorem \ref{theorem_evidence_for_conjecture}:

\begin{con}
Let $X'\subset \P(V)$ be an irreducible Legendrian subvariety.
Let $g \in \gotasp(V)$ be a nilpotent endomorphism  and $m$ be  an integer such that 
\mbox{$g^m\ne 0$} and $g^{m+1}=0$. 
Assume that the action of $\exp(tg)$ preserves $X'$.
Assume moreover, that $X'$ is singular at points of the image of the rational 
map $g^m(X')$. 
Then $X'$ is decomposable.
\end{con}

We also note the improved relation between projective automorphisms 
of a Legendrian subvariety and quadratic equations satisfied by its points:

\begin{cor}
Let $X\subset \P(V)$ be an irreducible Legendrian subvariety 
for which conjecture \ref{theorem_automorphisms_are_symplectic} holds
(for example $X$ is smooth).
If $G < \P\Gl(V)$ is the maximal subgroup preserving $X$, 
then $\dim G = \dim \I_2(X)$, where $\I_2(X)$ is the space of homogeneous quadratic 
polynomials vanishing on $X$.
\end{cor}
\begin{prf}
It follows immediately from the statement of the conjecture and \cite[lem. 5.6]{jabu}.
\end{prf}

Finally,
it is important to note, that theorem \ref{theorem_evidence_for_conjecture} part III.~does not imply 
that every torus acting on an indecomposable, but singular Legendrian variety  $X'$ 
is contained in the image of $\Sp(V)$. 
It only says that the intersection of such a torus with the weks-symplectic part is always finite. 
Therefore if there is a non-trivial torus acting on $X'$, 
there is also some non-trivial connected subgroup of $\Sp(V)$ acting on $X'$ 
and also some quadratic equations in the ideal of $X'$.

\section{Toric Legendrian Subvarieties}\label{section_toric}

Within this section $X$ is a toric subvariety of dimension $n-1$ in a projective space of dimension $2n-1$.
We assume it is embedded torically, so that the action of $T:=(\C^*)^{n-1}$ on $X$ 
extends to an action on the whole $\P^{2n-1}$,
but we do not assume that the embedding is projectively normal. 
The notation is based on \cite{sturmfels} though we also use technics of \cite{oda}.
We would like to understand when $X$ can be Legendrian with respect to some contact structure
on $\P^{2n-1}$ and in particular, when it can be a smooth toric Legendrian variety.

There are two reasons for considering non projectively normal toric varieties here. 
The first one is that the new example we find is not projectively normal. 
The second one is the conjecture \cite[conj. 2.9]{sturmfels}, 
which says that a smooth, toric, projectively normal variety is defined by quadrics. 
We do not expect to produce a counterexample to this conjecture 
and on the other hand all smooth Legendrian varieties defined by quadrics are known
(see \cite[thm.5.11]{jabu}).

\medskip

In addition we assume that either $X$ is smooth or at least the following condition is satisfied:
\begin{itemize}\label{star_condition}
\item[($\star$)]
The action of the torus $T$ on $\P^{2n-1}$ preserves the standard contact structure on $\P^{2n-1}$.
In other words, the image of $T\ra \P\Gl_{2n}$ is contained in the image of $\Sp_{2n} \ra \P\Gl_{2n}$.
\end{itemize}

In the case where $X$ is smooth, the ($\star$) condition is always satisfied
by theorem \ref{corollary_that_conjecture_is_true_for_smooth}. 
But for some statements below we do not need non-singularity, so we only assume ($\star$).

\begin{theo}\label{theorem_toric_legendrian}
Let $X\subset \P^{2n-1}$ be a toric (in the above sense) non degenerate Legendrian subvariety
satisfying ($\star$).
Then there exists a choice of symplectic coordinates\footnote{That is, coordinates for which the symplectic form has matrix 
$
\left(
\begin{array}{cc}
0 &\Id_n\\
-\Id_n & 0
\end{array}
\right)
$.}
and coprime integers $a_0 \ge a_1 \ge \ldots \ge a_{n-1} > 0$
such that $X$ is the closure of the image of the following map:
$$
T \ni (t_1, \ldots, t_{n-1}) 
\mapsto
[-a_0 t_1^{a_1} t_2^{a_2} \ldots t_{n-1}^{a_{n-1}}, \ \
a_1 t_1^{a_0}, \  a_2 t_2^{a_0},\  
\ldots,  \ a_{n-1}t_{n-1}^{a_0},$$
$$
 t_1^{-a_1} t_2^{-a_2} \ldots t_{n-1}^{-a_{n-1}}, \ \ 
 t_1^{-a_0}, \ t_2^{-a_0},\  
\ldots, \  t_{n-1}^{-a_0}] \in \P^{2n-1}.
$$
In other words, $X$ is the closure of the orbit of a point 
$$
[-a_0, a_1, a_2, \ldots a_{n-1}, 1, 1, \ldots 1]\in \P^{2n-1}
$$
under the torus action with weights 
$$
w_0:=(a_1, a_2,\ldots, a_{n-1}),
$$
$$   
w_1:=(a_0,0,\ldots, 0), \
w_2:=(0,a_0,0,\ldots, 0), \
\ldots, \
w_{n-1}:=(0,\ldots, 0, a_0) 
$$
$$
\textrm{and } 
-w_0, -w_1, \ldots, -w_{n-1}.
$$

Moreover every such $X$ is a non degenerate toric Legendrian subvariety.
\end{theo}

We are aware, that for many choices of the $a_i$'s from the theorem, the action of the torus on 
$X$ (and on $\P^{2n-1}$) is not faithful, 
so that for such examples a better choice of coordinates could be done. 
But we are willing to pay the price of taking a quotient of $T$ to get a uniform description.
One advantage of the description given in the theorem is that a part of it is almost 
independent of the choice of the $a_i$'s. It is the $(n-1)$ dimensional ``octahedron'' 
$\conv\{w_1,\ldots w_{n-1}, -w_1, \ldots -w_{n-1}\}\subset \Z^{n-1} \otimes \R$.

\begin{prf}
Assume $X$ is Legendrian with respect to a symplectic form $\omega$, that $X$ is non degenerate,
that the torus $T$ acts on $\P^{2n-1}$ preserving $X$ and satisfies ($\star$).
Replacing if necessary $T$ by some covering we may assume that $T\ra \P\Gl_{2n}$
factorises through a maximal torus $T_{\Sp_{2n}} \subset \Sp_{2n}$:
$$
T\ra T_{\Sp_{2n}} \subset \Sp_{2n} \ra \P\Gl_{2n}.
$$

This implies, that for an appropriate symplectic basis 
the variety $X$ is the closure of the image of the map $T\ra \P^{2n-1}$ given by:
$$
T \ni t \mapsto [x_0 t^{w_0},x_1 t^{w_1}\ldots, x_{n-1}t^{w_{n-1}}, t^{-w_0}, t^{-w_1}\ldots, t^{-w_{n-1}}] \in \P^{2n-1}
$$
where $x_i \in \C$, $w_i \in \Z^{n-1}$ and for
 $v= (v_1,\ldots v_{n-1})\in \Z^{n-1}$ we let \mbox{$t^v := t_1^{v_1} \ldots t_{n-1}^{v_{n-1}}$}.
 This means that $X$ is the closure of the $T$-orbit of the point%
\footnote{Note
that usually one assumes that this point is 
just [1,\ldots,1]. 
In our case we would have to consider non-symplectic coordinates. 
We prefer to deal with a more complicated point.
} 
$[x_0, \ldots x_{n-1}, 1, \ldots, 1]$
where $T$ acts with  weights 
$w_0, \ldots w_{n-1}, -w_0,\ldots, -w_{n-1}$.

Since $X$ is non degenerate, then the weights are pairwise different. 
Also the weights are not contained in any hyperplane in $\Z^{n-1}\otimes \R$,
because the dimension of $T$ is equal to the dimension of $X$ and we assume $X$ has an open orbit of
 the $T$-action. 
So there exists exactly one (up to scalar) linear relation:
$$
 - a_0 w_0 + a_1 w_1 + \ldots + a_{n-1} w_{n-1} = 0 .
$$
We assume that the $a_i$'s are coprime integers.
Permuting coordinates appropriately we can assume that $|a_0| \ge |a_1| \ge \ldots \ge |a_{n-1}| \ge 0$.
After a symplectic change of coordinates, we can assume without loss of generality that all the $a_i$'s are non negative by exchanging 
$w_i$ with $-w_i$ (and $x_i$ with $-\frac{1}{x_i}$) if necessary. 
Clearly not all the $a_i$'s are zero so in particular $a_0>0$ and hence
$$
w_0 = \frac{a_1 w_1 + \ldots + a_{n-1} w_{n-1}}{a_0}.
$$
Therefore, if we set $e_i:= \frac{w_i}{a_0}$ for $i \in \{1, \ldots, n-1 \}$, 
the points $e_i$ form a basis of a lattice $M$ containing all $w_i$'s. 
The lattice $M$ might be finer than the one generated by the $w_i$'s. 
Replacing again $T$ by a finite cover, 
we can assume that the action of $T$ is expressible in the terms of weights in $M$.
Then:
$$
\begin{array}{rcl}
w_0     & = & a_1 e_1 + \ldots + a_{n-1} e_{n-1},\\
w_1     & = & a_0 e_1\\
        & \vdots&\\
w_{n-1} & = & a_0 e_{n-1}
\end{array}
$$

It remains to prove three things: that $a_{n-1}>0$, 
that $x_i$'s might be chosen as in the statement of the theorem 
and finally that every such variety is actually Legendrian.
We will do all three together. 

The torus acts symplectically on the projective space, 
thus the tangent spaces to the affine cone are Lagrangian if and only if 
just one tangent space at a point of the open orbit is Lagrangian.
So take the point $[x_0, \ldots x_{n-1}, 1, \ldots, 1]$. 
The affine tangent space is spanned by the following vectors:
$$
\begin{array}{rrl}
v:=      &(x_0,\ \     x_1,\ \  x_2, \  \ldots, \ x_{n-1}, & 1,       1,    1,    \ldots,  1)\\
u_1:=    &(x_0a_1, \   x_1a_0,\ 0, \    \ldots, \quad 0,   &-a_1,     -a_0, 0,    \ldots,  0)\\
u_2:=    &(x_0a_2, \   0,  \ \  x_2a_0, \ldots, \quad 0,   &-a_2,      0    -a_0, \ldots,  0)\\
\vdots \\
u_{n-1}:=&(x_0a_{n-1} \ ,0, \ \ 0  \ \  \ldots, x_{n-1}a_0,&-a_{n-1},  0,    0,    \ldots, -a_0)\\
\end{array}
$$

Now the products are following:
$$
\omega(u_i, u_j) =0;  
$$
$$
\omega(u_i, v) = 2(x_0 a_i + x_i a_0).
$$
Therefore the linear space spanned by $v$ and $u_i$'s is Lagrangian if and only if:
$$
x_i = -x_0 \frac{a_i}{a_0}.
$$
In particular, since $x_i \ne 0$, the $a_i$ cannot be zero either.
After another conformal symplectic base change,
we can assume that $x_0= -a_0$ and then $x_i = a_i$.
 On the other hand, the above equation is satisfied for the variety in the
theorem. Hence the theorem is proved.
\end{prf}

Our next goal is to determine for which values of the $a_i$'s the variety $X$ is smooth. 
The curve case is not interesting at all and also very easy, so we start from $n=3$, 
i.e.~Legendrian surfaces.

\subsection{Smooth Toric Legendrian Surfaces}\label{section_toric_surfaces}

We are interested in knowing when the toric projective surface with weights of torus action
$$
w_0:=(a_1,a_2), \ w_1:=(a_0, 0), \ w_2:=(0, a_0),
$$
$$
 -w_0 = (-a_1, -a_2), \ -w_1=(-a_0,0), \ -w_2 = (0,-a_0) 
$$
  is smooth.
Our assumptions on the $a_i$'s are following:
\begin{equation}\label{inequalities_for_surface}
a_0\ge a_1\ge a_2 >0
\end{equation}
and $a_0, a_1,a_2$ are coprime integers.

\begin{figure}[htb]
\centering
\includegraphics[scale=0.75]{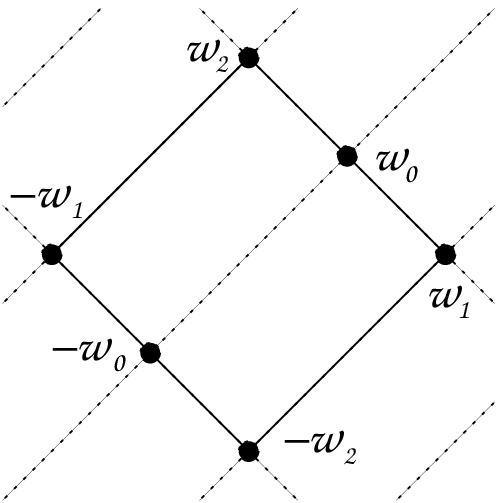}
\hspace{1cm}
\includegraphics[scale=0.75]{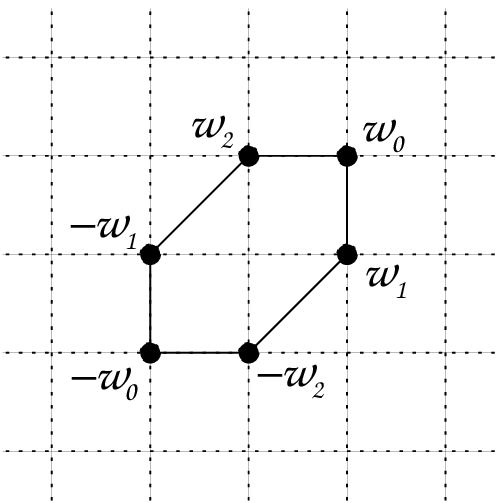}
\caption{
\footnotesize
The two examples of weights giving smooth toric Legendrian surfaces. 
}
 \label{figure_surface_examples}
\end{figure}

\begin{ex} \label{example_line_times_conic}
Let $a_0=2$ and  $a_1=a_2=1$ (see figure \ref{figure_surface_examples}). Then $X$ is the product of $\P^1$ and a quadric plane curve $Q_1$.
\end{ex}

\begin{ex}\label{example_plane_blown_up}
Let $a_0=a_1=a_2=1$ (see figure \ref{figure_surface_examples}).
Although the embedding is not projectively normal (we lack the weight $(0,0)$ in the middle),
the image is smooth anyway.
Then $X$ is the blow up of $\P^2$ in three non-colinear points.
\end{ex}

We will prove there is no other smooth example.

\medskip

\begin{figure}[htb]
\centering
\includegraphics[scale=0.5]{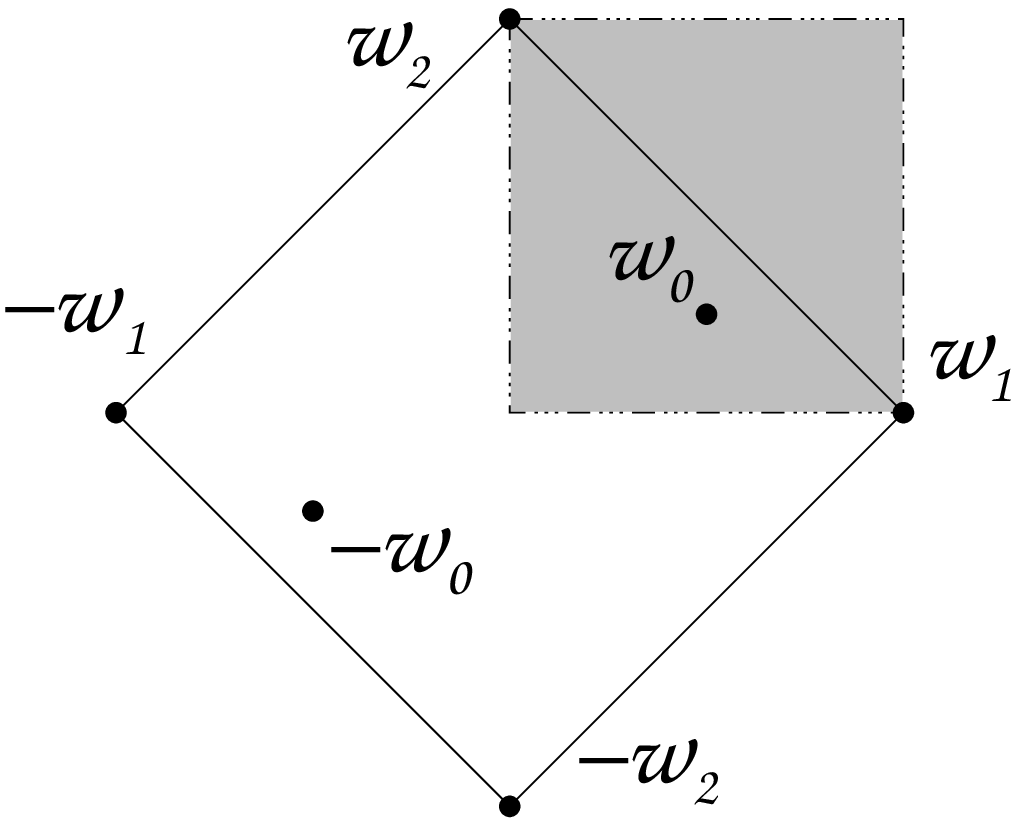}
\hspace{1cm}
\includegraphics[scale=0.5]{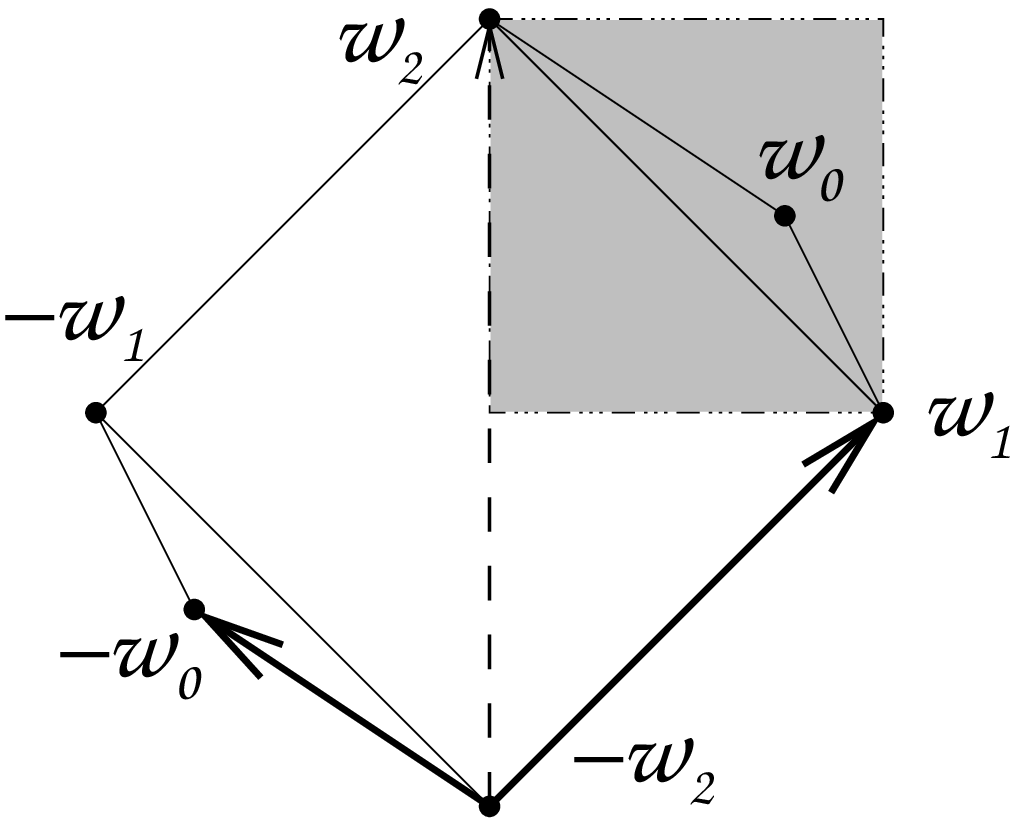}
\caption{
\footnotesize
Due to the inequalities $a_0 \ge a_1> 0$  and $a_0 \ge a_2> 0$,  
the weight $w_0$ is located somewhere in the gray square. 
The two cases we consider are if $w_0$ is also inside the square $\conv \{ w_1,w_2, -w_1, -w_2 \}$ 
(left figure) or it is outside (right figure). 
In the second case, a necessary condition to get a smooth variety, is that the two bold vectors generate a lattice containing all the weights. 
In particular the dashed vector can be obtained as an integer combination of the bold ones.
}
 \label{figure_surface_proof}
\end{figure}

We must consider two cases (see figure \ref{figure_surface_proof}):
either $a_0 > a_1 +a_2 $ (which means that $w_0$ is in the interior of the square 
$\conv \{ w_1,w_2, -w_1, -w_2 \}$) or 
$a_0 \le a_1+a_2$ (so that $w_0$ is outside or on the border of the square).

\smallskip

Geometrically, case $a_0 > a_1 +a_2 $ means, that the normalisation of $X$ is $\P^1\times \P^1$. 
It is just an easy explicit verification 
that $X$ is not smooth with these additional weights in the interior.

\smallskip

In the other case, 
for a vertex $v$ of the polytope 
$$
\conv \{ w_0,w_1,w_2, -w_0,-w_1, -w_2 \}
$$ 
we define the sublattice $M_v$ to have the origin at $v$ and to be generated by
$$
\{ w_0-v,w_1-v,w_2-v, -w_0-v,-w_1-v, -w_2-v\} .
$$
Since $X$ is smooth, for every vertex $v$
the vectors of the edges meeting at $v$ 
must form a basis of $M_v$ 
(compare with \cite[prop.2.4 \& lemma 2.2]{sturmfels}). 
In particular, if $v=-w_2$ 
(it is immediate from inequalities (\ref{inequalities_for_surface}) that $v$ is indeed a vertex), 
then $w_2-(-w_2) = (0, 2 a_0)$ can be expressed as an integer combination of 
$w_1+w_2= (a_0,a_0)$ and $-w_0+w_2 = (-a_1, a_0-a_2)$ 
(see the righthand side of figure \ref{figure_surface_proof}).
So write:
\begin{equation}\label{equation_for_surface_case}
(0, 2 a_0) = k (a_0,a_0) + l (-a_1, a_0-a_2)
\end{equation}
for some integers $k$ and $l$.
 It is obvious that $k$ and $l$ must be strictly positive, 
since $w_2$ is in the cone generated by $w_1+w_2$ and  $-w_0+w_2$ with the vertex at $-w_2$.
But then (since $a_0-a_2 \ge 0$) 
from equation (\ref{equation_for_surface_case}) on the second coordinate we get that either
 $k=1$ or $k=2$.

\smallskip

If $k=1$, then we easily get that:
$$
\left\{
\begin{array}{l}
a_0 = l a_1 \\
a_0 = a_1 +a_2
\end{array}
\right.
$$
Hence $(l-1)a_1 = a_2$ and by inequalities (\ref{inequalities_for_surface}) we get $l=2$ 
and therefore (since the $a_i$'s are coprime)
$(a_0,a_1, a_2) =(2,1,1)$ which is  example \ref{example_line_times_conic}.

\smallskip

If on the other hand $k=2$, then 
$$
a_0 = a_2 
$$
and hence by inequalities (\ref{inequalities_for_surface}) and since the $a_i$'s are coprime,
 we get $(a_0,a_1, a_2) =(1,1,1)$, which is example \ref{example_plane_blown_up}.

\begin{cor}\label{smooth_toric_surfaces}
If $X\subset \P^5$ is smooth toric Legendrian surface, 
then it is either $\P^1 \times Q_1$ 
or $\P^2$ blown up in three non-colinear points or plane $\P^2\subset \P^5$.
\end{cor}
\noprf


\subsection{Higher dimensional toric Legendrian varieties}
\label{section_higher_dimensional}

In this subsection we assume that $n\ge 4$. 
By means of the geometry of  convex bodies we will prove there is only one smooth 
toric non degenerate Legendrian variety   
in dimension $n-1=3$ and no more in higher dimensions.
We use theorem \ref{theorem_toric_legendrian} so that we have a toric variety with weights:
$$
\begin{array}{l}
w_0:=(a_1,a_2,\ldots, a_{n-1}),\\
w_1:=(a_0,0,\ldots 0),\\
\vdots\\
w_{n-1}:=(0,\ldots 0, a_0),\\
-w_0, -w_1, \ldots, -w_{n-1}
\end{array}
$$ 
where the $a_i$'s are coprime positive integers with $a_0\ge a_1 \ge \ldots \ge a_{n-1}$.

\begin{figure}[htb]
\centering
\includegraphics[scale=0.7]{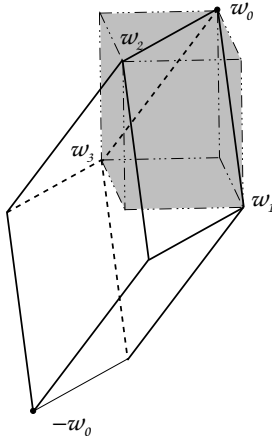}
\caption{
\footnotesize
The smooth example in dimension 3: $(a_0,a_1,a_2,a_3)= (1,1,1,1)$.
}
 \label{octahedron_case4}
\end{figure}

\begin{ex}\label{example_line_cubed}
Let $n=4$ and $(a_0,a_1,a_2,a_3)= (1,1,1,1)$. 
Then the related toric variety is $\P^1 \times \P^1 \times \P^1$ 
(see figure \ref{octahedron_case4}).
\end{ex}

Further, let $A$ be the polytope defined by the weights:
$$
A:=\conv \{ w_0,w_1,\ldots , w_{n-1}, -w_0, -w_1, \ldots, -w_{n-1}\} \subset \Z^{n-1} \otimes \R.
$$

\begin{lem}\label{edges_of_A}
Let $I,J \subset \{1, \ldots, n-1\}$ be two complementary subsets of indexes. 
\begin{itemize}
\item[(a)]
Assume $i_1, i_2 \in I$ and $ i_1 \ne i_2$. If 
$$
\left| \sum_{i \in I} a_i -\sum_{j \in J} a_j \right| < a_0 ,
$$
then the interval $(w_{i_1},w_{i_2})$ is an edge of $A$.
\item[(b)]
Assume $k \in I$ and $l \in J$. If
$$
\sum_{i \in I} a_i -\sum_{j \in J} a_j  > a_0, 
$$
then both intervals $(w_0,w_k)$ and $(w_0, -w_l)$ are edges of $A$.
\item[(c)]
If $k,l \in \{1, \ldots, n-1\}$ and $k \ne l$, then  $(w_k, -w_l)$ is an edge of $A$.
\end{itemize} 
\end{lem}

\begin{prf}
Fix $\epsilon >0$ small enough,
set $\alpha:= \sum_{i \in I} a_i -\sum_{j \in J} a_j$
and define the following hyperplanes in $\Z^{n-1} \otimes \R $:
$$ 
H_a:= \left\{ 
\sum_{i \in I} x_i - (1-\epsilon)  \sum_{j \in J} x_j  = a_0 
\right\},
$$
$$
H_b:=
\left\{
\left(a_0-a_k\right)\left(\sum_{i \in I} x_i - \sum_{j \in J} x_j -\alpha \right)
+(\alpha -a_0)\left(x_k-a_k\right) =  0
\right\},
$$
$$
H'_b:=
\left\{
\left(a_0+a_l\right)\left(\sum_{i \in I}  x_i - \sum_{j \in J} x_j - \alpha  \right)
+(\alpha -a_0)\left(x_l+a_l\right) =  0
\right\}
$$
$$
\textrm{and }
H_c:=
\left\{
x_k-x_l =  a_0
\right\}.
$$

Assuming the inequality of (a), 
$H_a \cap A$ is equal to $\conv\{w_i \mid i \in I \}$ and 
the rest of $A$ lays on one side of $H_a$.
So $H_a$ is a supporting hyperplane for the face \mbox{$\conv\{w_i \mid i \in I \}$},
which is a simplex of dimension $(\#I-1)$ 
and therefore all its edges are also edges of $A$ as claimed in (a).

\smallskip

Next assume that the inequality of (b) holds. 
Then $H_b$ (respectively $H'_b$) is a supporting hyperplane for the edge $(w_0,w_k)$
(respectively $(w_0,-w_j)$).  

\smallskip

Similarly, in the case of (c), $H_c$ is a supporting hyperplane for $\{w_k, -w_l\}$.

\end{prf}

\begin{theo}\label{classification_toric_higher_dimensional}
Let $X\subset \P^{2n-1}$ be a toric non degenerate Legendrian variety of dimension $n-1$
satisfying ($\star$) (see page \pageref{star_condition}).
If $n\ge 4$ and normalisation of $X$ has at most quotient singularities,
 then $n=4$ and $X=\P^1 \times \P^1 \times \P^1$.
\end{theo}

\begin{prf}
Since the normalisation of $X$ has at most quotient singularities,
it follows that the polytope $A$ is simple, 
i.e.~every vertex has exactly $n-1$ edges (see \cite{fulton} or \cite[\S2.4, p. 102]{oda}).
We will prove this is impossible, unless $n=4$ and $(a_0, a_1, a_2,a_3)= (1,1,1,1)$.

\smallskip

If  $w_0 \in B:=\conv \{w_1,\ldots , w_{n-1}, -w_1,\ldots -w_{n-1}\}$,
then $A$ is just equal to $B$ 
 and clearly in such a case every vertex of $A$ has  $2(n-2)$ edges. 
Hence more than $n-1$ for $n\ge 4$.

\smallskip
 
Hence from now on we can assume that $a_1 + \ldots + a_{n-1} > a_0$. 
So by lemma \ref{edges_of_A}(b), $(w_0, w_i)$ is an edge for every $i \in \{1, \ldots, n-1\}$.

Choose any $j \in \{1, \ldots, n-1\}$ and set $I:=\{1, \ldots,j-1, j+1, \ldots, n-1\}$.

If either  
$$
\left|\left(\sum_{i \in I} a_i\right) \ -  \ a_j \right| \ < \ a_0 \quad \textrm{ or}
$$
$$
\left(\sum_{i \in I} a_i\right) \ -  \ a_j \  > \ a_0, 
$$
then using lemma \ref{edges_of_A} we can count the edges at either $w_i$ or $w_0$ 
and see that there is always more than $n-1$ of them. 
We note that $a_j -\left(\sum_{i \in I} a_i\right) \ge a_0$ never happens due to our assumptions on the $a_i$'s.

\smallskip

Therefore the remaining case to consider is
$$
\left(\sum_{i \in I} a_i\right) \  - \  a_j \   = \  a_0, 
$$
where the equality holds for every $j \in \{1, \ldots,n-1\}$.
This implies:
$$
a_1 =a_2 =\ldots = a_{n-1} = \frac{1}{n-3} a_0.
$$
Since the $a_i$'s are positive integers and coprime, we must have 
$$
(a_0,a_1, \ldots a_{n-1}) = (n-3, 1,\ldots, 1)
$$ 
which is exactly example \ref{example_line_cubed} for $n=4$.
Otherwise, if  $n\ge 5$ we can take $J:=\{j_1,j_2\}$ for any two different 
$j_1,j_2 \in \{1, \ldots, n-1\}$ and set $I$ to be the complement of $J$.
Then $\# I \ge 2$ and  by lemma \ref{edges_of_A}(a) and (c) there are  too many of edges 
at the $w_i$'s.
\end{prf}

\begin{cor}\label{smooth_toric_varieties}
If $X\subset \P^{2n-1}$ is a smooth toric Legendrian subvariety 
and $n\ge 4$,
then it is either a linear subspace or $n=4$ and $X= \P^1 \times \P^1 \times \P^1$.
\end{cor}
\noprf


\end{document}